\documentclass[11pt,a4paper]{article}
\usepackage[T1]{fontenc}
\usepackage[utf8]{inputenc}
\usepackage[margin=27mm]{geometry}
\usepackage{amsmath,amssymb,amsthm,mathtools,mathrsfs,booktabs,microtype}
\usepackage[colorlinks=true,linkcolor=blue,citecolor=blue,urlcolor=blue]{hyperref}
\hypersetup{pdftitle={An exhaustion-aware product moment and zero-free half-planes for Hecke L-functions},pdfauthor={Xinrui Chang, HaoMing Wang}}
\allowdisplaybreaks
\numberwithin{equation}{section}
\newtheorem{theorem}{Theorem}[section]
\newtheorem{proposition}[theorem]{Proposition}
\newtheorem{lemma}[theorem]{Lemma}

\theoremstyle{definition}

\newcommand{\eps}{\varepsilon}
\newcommand{\Nm}{\mathrm N}
\newcommand{\rad}{\operatorname{rad}}

\newcommand{\Rea}{\operatorname{Re}}
\newcommand{\calS}{\mathcal S}
\title{An exhaustion-aware product moment and zero-free half-planes for Hecke $L$-functions}
\author{Xinrui Chang\\\small Fudan University\\\small\texttt{25800180002@m.fudan.edu.cn}
\and HaoMing Wang\\\small Fudan University\\\small\texttt{25800180121@m.fudan.edu.cn}}
\date{}
\begin{document}
\maketitle
\begin{abstract}
We refine the ordinary product moment in the cubic-theta comparison by allowing a returned prime product to be exhausted. The return stops either when the positive-prime length inequality holds or when the established prime-free ordinary estimate applies. Keeping the extracted central coefficients gives a local normalization inequality on a larger domain. On the critical ordinary-length interval this yields a net first-mean gain of $1/20000$. We prove the complementary-class, short-sum and algebraic comparison inequalities directly, and obtain the endpoint $B_*-1/160000=0.87446205769353\ldots$, where $B_*$ is Liu's algebraic boundary. The zero-free conclusion is conditional on the cited general-field theta, detector, prior-bound and Mellin input theorems in their stated uniform forms. The proof includes the actual return weights, both approximate-functional-equation arms, the principal-residue normalization and descent to arbitrary number fields.
\end{abstract}

\section{Introduction}
The cubic-theta method relates zero-free half-planes to moments of ordinary character sums and to the number of character parameters in a detected amplitude class. Liu's general-field formulation~\cite{Liu} supplies the ordinary moment and the zero-detecting comparison. His subsequent amplified comparison~\cite{H} reaches the algebraic boundary $B_*=0.87446830769353045955\ldots$. We enlarge the admissible additional prime product by changing the stopping rule in the ordinary return. The change concerns the complete positive norm, before restriction to a detected class.

Let $r_*$ be the unique root in $(9/20,23/50)$ of
\begin{equation}\label{eq:cubic}
 p(r)=218582r^3-376469r^2-11376r+61560=0.
\end{equation}
Define
\begin{equation}\label{eq:root}
 \ell_*={3410r_*^2-7127r_*+2520\over2286r_*^2-5011r_*+1800},\qquad
 B_*={11\over12}-{\ell_*\over4},\qquad k_*={5\over6}-{\ell_*\over2}.
\end{equation}
For independent variables $r,k$, put
\begin{align*}
 n_b(r,k)&=(-6858r^2+15033r-5400)k+600r^2-1837r+720,\\
 n_e(r,k)&=(-702r^2+34995r-40500)k-408r^2-3775r+5400,\\
 n_0(r,k)&=(4806r^2-10887r+8100)k-696r^2+1403r-1080,
\end{align*}
and set
\begin{equation}\label{eq:b}
 b_*=-{\ell_*\partial_r n_e(r_*,k_*)+\partial_r n_0(r_*,k_*)\over
                   \partial_r n_b(r_*,k_*)}.
\end{equation}
The derivatives in this definition hold $k$ fixed. Our new length and boundary are
\begin{equation}\label{eq:boundary}
 \ell=\ell_*+{1\over40000},\qquad B={11\over12}-{\ell\over4}
       =B_*-{1\over160000}.
\end{equation}

\begin{theorem}\label{thm:main}
Assume the general-field analytic inputs specified in Section~\ref{sec:inputs}. For each fixed number field $F$ containing $\mu_{12}$, every finite-order Hecke $L$-function over $F$ is zero-free in $\Rea s>B$. If these inputs hold for every such $F$, the same conclusion holds for every finite-order Hecke $L$-function over every number field. In particular, it holds for the Riemann zeta function and for Dirichlet $L$-functions under those inputs.
\end{theorem}

The main new estimate is Theorem~\ref{thm:norm}. Its proof retains the signed divisor expansion until the actual ordinary children have been formed. The pointwise prime cap is then applied to whole factors. At exhaustion the remaining product is estimated by the prime-free theorem, even when its formal positive-prime defect is still positive. An exact normalization identity pays this operation. The numerical improvement follows by using a slightly longer selected product on $41/100\le m\le43/100$.

The terminal ordinary estimates and theta construction are invoked from the cited papers. The local return, the scalar comparison, and the passage from a positive complete saving to the new endpoint are developed below. In particular, the scalar inequalities are proved by rational differentiation and a double-root identity rather than by an auxiliary certificate table. All constants, derivative orders and thresholds may depend on the fixed field and target character; the boundary $B$ does not.

\section{Arithmetic families and analytic inputs}\label{sec:inputs}
Fix $F\supset\mu_{12}$. Choose a finite set $S$ of places, admissible generators of the relevant ideals, and a finite group $\Theta$ of finite-place and unit characters as in~\cite[Part III]{Liu}. Write $q_a=\Nm a$. The sextic residue character $\chi_a$ is extended by zero at every nonunit of its displayed modulus. The displayed moving radical is retained throughout a divisor extraction, including factors that cancel algebraically. We use the normalized full Gauss sum
\[
 \mathcal G(a,h)=q_a^{-1/2}\sum_{v\bmod a}\chi_a(v)e(hv/a),\qquad
 \mathcal G(1,h)=1,
\]
with the trace additive character and local measures of~\cite{Liu}. The imprimitive formula is used also at nonunit frequencies.

For a primitive periodic character $\lambda$ of modulus $\mathfrak m$ and $(D,\mathfrak m)=1$, the nonzero part of orbit Poisson, with physical norm scale $Y$ and $\Phi(0)=0$, is
\begin{equation}\label{eq:poisson}
 \sum_{d\mid D}{\mu(d)\lambda(d)Y\tau(\lambda)\over q_dq_{\mathfrak m}}
       \sum_{[h]\ne0}\overline{\lambda(h)}
       K_\Phi\left({Yq_h\over q_dq_{\mathfrak m}}\right).
\end{equation}
This follows by expanding coprimality, writing the row as $dk$, applying lattice Poisson to each residue class, and evaluating its finite Fourier transform. If the primitive modulus is one, add the origin
$Y\operatorname{vol}(G/\Gamma)\widehat\Phi(0)\sum_{d\mid D}\mu(d)/q_d$; otherwise that origin vanishes. For an imprimitive character retain its full periodic Fourier transform instead of replacing it by $\tau(\lambda)\overline\lambda$. The origin has the compact-quotient measure, whereas the nonzero orbits unfold to $G$. These are the two measures in the theta normalization of~\cite{Liu}.

An ordinary annular factor has the form
\begin{equation}\label{eq:ordinary}
 S_{X,V}(u)=X^{-1/2}\sum_a c(a)\chi_u(a)V(q_a/X),
\end{equation}
where the sum is over the prescribed generator classes, $c$ belongs to the fixed finite character span, and the common coprimality conditions are written with the same zero extension. Different factors may have different admissible profiles and independent Mellin heights. A prime factor is the corresponding normalized sum on a full prime annulus, with the original finite-$\Theta$ coefficients. Its logarithmic width relative to $U$ is denoted by $z_i$. Products below always use disjoint prime slots. These definitions do not allow coefficients chosen as a function of the current row.

We recall precisely the external results used. The theta coefficients $a_0$, their twisted multiplicativity, reciprocity, full local Gauss identities and orbit Poisson summation are those of~\cite[Part III]{Liu}, with the metaplectic normalization originating in~\cite{KP,KPcorr}. The inverse moment and the ordinary product moment, including arbitrary bounded ordinary lengths when no prime factor remains, are the theorems in the sections ``The unweighted product moment'' and its finite-width induction in~\cite{Liu}. They include fixed divisor weights, natural zeros, finite twists and polynomial control of all requested finite-order profile seminorms.

For zero detection we use the original inverse detector, its amplitude-class count, the same-object contour comparison, both actual-conductor approximate functional equations, and the direct estimate of~\cite{Liu,H}. The prior result $\beta_*(F)\le B_*<0.8745$ from~\cite{H} is assumed before any new prime cap is used. Here
\[
 \beta_*(F)=\sup\{\Rea\rho:L_F(\rho,\eta)=0,\ 0<\Rea\rho<1,
                    \ \eta\text{ of finite order}\}.
\]
This prior bound licenses the local cap with $\kappa=749/1000$ and the count with parameter $k_*$. Its proof is an input, rather than an inference from the improved endpoint.

Finally we use the normalized theta/Mellin identity in~\cite[the principal contour and continuation argument]{Liu,H}. With $c_0=(4+b)/6$ it identifies the Mellin transform of $f_\eta(Z)$ with
\begin{equation}\label{eq:mellin}
 \exp((s-5/6)^2){H_\eta(s)\over L_F^S(s,\eta)},
\end{equation}
where $H_\eta$ is holomorphic and nonzero for $\Rea s>437/500$. The direct estimate for the normalized average $J_\eta$ has exponent $B(\ell)-c_0+\eps$, and the principal contour produces the same $f_\eta$. These input statements include their local corrections and their allowed order of choices. Theorem~\ref{thm:main} is an implication from these analytic results; it does not require the reader to adopt unrecorded numerical computations.

\section{Return weights and profile closure}\label{sec:return}
We first describe why the two transforms return the families just defined. This also fixes how their weights are used in the ordinary argument.

\subsection{The inverse family}
The canonical inverse energy of~\cite{Liu} is
\begin{equation}\label{eq:inverse}
 \mathcal E=Z^{-N-V}\sum_{\substack{f\text{ squarefree}\\q_f\asymp Z^V}}w(f)
 \sum_{0<q_k\le CZ^M}\left|\sum_{n\text{ squarefree}}
 a_0(n)\nu(n)\rho(n)\chi_n(k)\chi_n(f)^4W(q_n/Z^N)\right|^2,
\end{equation}
where $\rho(n)=1_{(n,R_\rho)=1}$, $w(f)\ll\tau(f)^L$ is independent of $k$, and $\nu$ is a fixed finite twist. Complete cube extraction writes $N=r+3a$ and gives the original squared normalization $Z^{-N-V}Y_b$, with $Y_rY_b^3=Z^N$. Möbius coefficients in this extraction remain signed.

In the first Poisson transform write $n_i=A_iCu_i$ and expand $1_{(u_1,u_2)=1}$ by its full divisor identity. The cube-prime valuation states give ideals $q_0,J=sJ_2,E,R_1$, with $b_1b_2=q_0^2J_2^2s$. The actual residual coefficient is
\[
 \mu(u)\nu_1(u)\rho(u)\overline{\chi_u(y)}\chi_u(CJ)^4\chi_u(d_k)
                 1_{(u,q_0)=1},\qquad y=hf^2E.
\]
For fixed $y$ the old labels $f,h$ have divisor-bounded multiplicity, since $f^2E\mid y$ determines $h$. This bound is taken on the original support, before positive enlargement. Thus it does not incur a second factor $Z^V$.

For the second transform write $x_j=g'z_j$, $z_j=v'n_j$, $d_2\mid g'$, and put
\[
 k'=d_kd_2k'',\qquad f'=JCd_2v',\qquad
 \rho_\gamma(n)=\rho(n)1_{(n,\rad(q_0t'g'/d_2))=1}.
\]
The returned coefficient is exactly
\begin{align*}
 &a_0(n)\nu_2(n)\rho_\gamma(n)\chi_n(d_kk'')
            \overline{\chi_n(d_2)}\chi_n(CJv')^4\\
 &\hspace{35mm}=a_0(n)\nu_2(n)\rho_\gamma(n)\chi_n(k')\chi_n(f')^4.
\end{align*}
Indeed $\chi_n(d_2)^5=\overline{\chi_n(d_2)}$ when the arguments are coprime; otherwise both sides vanish. The fixed unit Gauss factors give $\nu_2=\nu\overline{\xi_1}\xi_2$. On the support forced by twisted multiplicativity, $J,C,d_2,v'$ are pairwise coprime and squarefree. Hence $f'$ is an admissible label, including at every natural zero.

\begin{lemma}\label{lem:projection}
The return uses one common positive structural projection measure. If each child satisfies $\mathcal E_{\rm child}\le A Z^{F'+\lambda}\calS_{\rm child}$ uniformly over the permitted data, the parent bound contains $A Z^\lambda$ once.
\end{lemma}
\begin{proof}
Fix $\gamma=(q_0,t',g'/d_2)$. For fixed $(\gamma,f',k')$, splitting $f'$ into $J,C,d_2,v'$ and then splitting $J$ has divisor-bounded multiplicity; the old $(f,h)$ fiber was already counted. If $j,t,g,\vartheta,R$ denote the respective logarithmic scales of $J,t',g',d_2,R_1$, the structural mass exponent is
\[
 C_c=a+R/2-j+t+g-\vartheta.
\]
Write its weights as $A_\gamma$ and use the same measure
\[
 d\nu(\gamma)=Z^{-C_c}\sum_\gamma A_\gamma\delta_\gamma
\]
in both Cauchy integrals. The column and label centers are
\[
 N_c=r-A_i-B-t-g-v,\qquad V_c=B+\vartheta+v+j,\quad F_c=N_c+V_c.
\]
The exact Poisson prefactors satisfy
\begin{equation}\label{eq:normalization}
 \kappa_i+\lambda_c+C_c+2F_c=r+3a+V=F.
\end{equation}
One $F'$ restores the normalization in~\eqref{eq:inverse}; the other comes from the child bound. If the two resulting positive integrals are $I_i=\int\mathcal E_i\,d\nu_i$ and $T_i=\int\calS_i\,d\nu_i$, use
\[
 (I_1I_2)^{1/2}\le A Z^{F'+\lambda}(T_1T_2)^{1/2}.
\]
The measures and twists on the two sides need not coincide. Their common uniform bound and~\eqref{eq:normalization} transmit $A$ and $\lambda$ once. In particular, no interchange of the geometric mean with the integrals is used.
\end{proof}

\subsection{Profiles and the ordinary return}
All differentiation is performed in actual fixed-ratio coordinates $u_i=\log(q_i/Y_i)$. A typical separated profile is
\[
 P_Z(u)=\chi(u)K(Z^\zeta e^{L(u)})e^{R(u)}
             \prod_jW_j(c_je^{L_j(u)})e^{i\mathbf t\cdot Au}.
\]
The norm-sheet kernels from the fixed Schwartz--Bruhat lifts obey logarithmic bounds as their argument tends to zero and arbitrary polynomial decay at infinity. Integration by parts on the fixed $u$-box therefore gives, for each prescribed $A$, some finite $J,C$ such that
\begin{equation}\label{eq:profile}
 \int(1+|\xi|)^A|\widehat P_Z(\xi)|\,d\xi
 \ll(\log Z)^{d_F}(1+|\mathbf t|)^C\calS_J.
\end{equation}
The bound is uniform in the actual divisor scales $c_j$. It has no factor $Z^{J\eta}$ from a coarse logarithmic mesh.

The physical positive row lift need not be compactly supported. Apply the second exact Poisson transform to that lift before identifying a child ball. Its compact dual support gives the actual upper bound on the new row. A negative formal center alone does not delete a nonempty bounded row.

We spell out the ordinary paired weights. In a final artificial-divisor dyad $U^{t_-}\le q_t<2U^{t_-}$, let $\widetilde r_j$ be the logarithmic removed ordinary length, $\omega_j$ the fixed removed prime-slot length, $e_j$ the support-conversion error, and $\pi_j,\pi_{0j}$ the nominal endpoint allowances on side $j$. The actual central multiplier on that side is
\[
 U^{(-r_j^++e_j+\pi_{0j})/2},\qquad
 r_j^+=\widetilde r_j-(\pi_j-\pi_{0j}),\quad
 \widetilde r_j+\omega_j\ge t_-.
\]
Thus if $c_{t,\alpha}$ is the frozen coefficient of the complete last-divisor expansion, the weight in both positive integrals is
\[
 b_{t,\alpha}=|c_{t,\alpha}|
 U^{(-r_1^++e_1+\pi_{01}-r_2^++e_2+\pi_{02})/2}.
\]
The original endpoint rules give
\[
 t_-+\log_U(b_{t,\alpha}/|c_{t,\alpha}|)
 \le{e_1+\omega_1+e_2+\omega_2+\pi_1+\pi_2\over2}\le2\theta_N.
\]
Here $\theta_N=O(1/\log U)$ records the total fixed-ratio endpoint width of the finite factor list. With divisor multiplicities included, $\sum b_{t,\alpha}\ll U^{2\theta_N+\eps}$. The common structural radical $sB$ contributes its actual count $U^{p_2-s_0+\eps}$ and its opposite normalization $U^{-(p_2-s_0)}$. Hence each Cauchy-side measure is the structural measure times this same $b_{t,\alpha}$ and the absolute Fourier density. Integrate that density using~\eqref{eq:profile}, then insert the same low-rank uniform estimate separately on the two sides. In particular, the factor $U^{t_-}$ is never charged again as a bare divisor count.

Reflection also respects row dependence. For primitive inducing character $\psi_k^*$, redundant radical $R_{0,k}$ and actual conductor $C_k$, the Euler expansion of the functional equation is
\begin{align*}
 S_k^0(X;W)=\epsilon_k
 \sum_{\substack{d\mid R_{0,k}\\\rad(h)\mid R_{0,k}}}
 {\mu(d)\psi_k^*(d)\overline{\psi_k^*(h)}\over\sqrt{q_dq_h}}
 S_{\overline{\psi_k^0}}^0\left({C_kq_d\over Xq_h};W^\sharp\right).
\end{align*}
Its absolute Euler mass is bounded, separately on each row, by
\[
 L_k\le\prod_{p\mid R_{0,k}}{1+q_p^{-1/2}\over1-q_p^{-1/2}}
       \ll_\eps U^\eps.
\]
To see uniformity, first fix a sufficiently large prime cutoff; above it the logarithm of each factor is at most $\eps\log q_p$, and below it only a fixed constant remains. Decompose $W^\sharp$ into actual annuli. Their central coefficients are summable, including the logarithmic small-argument terms, and depend on the incoming profiles rather than the row. For each profile apply rowwise Cauchy with $L_k$ before the common two-scale Sobolev estimate. This retains the mixed derivative and does not require a common Euler-divisor set for all rows. Equation~\eqref{eq:profile} now closes any prescribed finite-order child seminorm. The strict width decrease bounds the depth in advance, so these orders can all be fixed before the external height exponent.

\section{An exhaustion-aware ordinary norm}\label{sec:norm}
Put
\[
 \kappa={749\over1000},\qquad K=6\kappa={2247\over500},\qquad c_K={K-1\over K}.
\]
The prior half-plane in Section~\ref{sec:inputs} permits the pointwise squared prime cost $U^{\kappa z+\eps}$. For $m\in[m_-,r_o]$ form
\[
 P=S_{U^m,V_1}S_{U^m,V_2}Q_0,\qquad P^c=P-P^{\rm cmp}.
\]
The comparison replaces the ordinary scales $(m,m)$ by $(1/4,2m-1/4)$, preserving their product, both profiles, and all prime factors. The positive row measure is the prescribed smooth measure $W_{\rm row}(q_u/U)$ of~\cite{H}.

\begin{theorem}\label{thm:norm}
Assume the local arithmetic and ordinary input theorems of Section~\ref{sec:inputs}. Suppose
\begin{equation}\label{eq:domain}
 {41\over100}\le m_-\le r_o\le{43\over100},\quad 0\le W\le{1\over40},\quad
 1-2r_o-7W>0,\quad {5\over2}-4m_-+KW<1.
\end{equation}
Let $Q_0,E$ be complete disjoint prime products of widths
$0\le z\le(1-2m)/K$ and $0\le w\le W$. Then, for every $\eps>0$,
\begin{align}\label{eq:norm}
 \sum_u W_{\rm row}(q_u/U)|P^c(u)E(u)|^2&\ll U^{1+\eps}\calS_J(1+H)^C,\\
 \sum_{u\in\mathcal U} W_{\rm row}(q_u/U)|P(u)E(u)|^2&\ll U^{1+\eps}\calS_J(1+H)^C.\nonumber
\end{align}
Here $\mathcal U$ is the source sixth-power-free nonunit family with primitive inducing character outside $\Theta$. The estimates retain the independent profiles, finite twists, natural zeros and common punctures. They are uniform on fixed compact parameter families with positive margins in~\eqref{eq:domain}.
\end{theorem}

\subsection{Valuation allocation and completed lengths}
Expand the complete signed coefficient product first. The valuation classes in the full Gauss correlation formulas~\cite[the Gauss-norm construction]{H} partition the added width as
\[
 w=w_f+f_0+f_1+u,\quad u=g_{\rm good}+v+c,\quad t=f_0+u,
 \quad L_1=f_1+4T_1+e_1,\quad L=v+4T+e,
\]
where all variables are nonnegative. Write $A_t=2m+z+w$. The completed lengths and coefficient exponent are
\begin{align}\label{eq:completed}
 A^*&=A_t-f_1-t-L_1-L-5c,\\
 M^*&=1-6t-f_1-L_1-L-5c+T_1+5v+T,\nonumber\\
 q^*&=f_1+v+c,\qquad m^*=M^*-q^*,\nonumber\\
 R_{\rm loc}&=-f_0-L_1+T_1-S_1-u-L+v-4c+T-S.\nonumber
\end{align}
For the original common-support ideals $C,D$, put
$\mu_{\rm com}=\min(\log_Uq_C,\log_Uq_D)$, so that $w_f\le\mu_{\rm com}$.
Let $z_C$ be the extracted part of the initial product and
$\eta_K=1-2m-Kz\ge0$. The ordinary positive-prime defect is exactly
\begin{equation}\label{eq:defect}
 D_{\rm free}=A^*+(K-1)(z-z_C)-M^*
 =6w-\eta_K-(K-1)z_C-5w_f-5f_1-T_1-5v-T.
\end{equation}
No length inequality has been used in this identity. Ordinary removal gives $L_1+L+5c\le2m+o(1)$, whence
\[
 M^*\ge1-2r_o-6W-o(1),\qquad m^*\ge1-2r_o-7W-o(1)>0.
\]
Support edges are chosen inside this fixed margin.

Keep the actual outer ideal $D$ and its Fourier measure. The numerical benchmark $\widehat d=\min(\log_Uq_D,A^*)$ changes its complementary exponent by the opposite amount; it does not introduce a second family of $D$'s. With the source complementary divisor exponents $R,E_0$, split the actual formal ordinary length
\[
 m_f=m^*-\log_Uq_C+\widehat d-R-E_0
\]
at $1/100$. The large-formal branch has an actual ordinary return. The other branch will be bounded directly.

\subsection{Removing whole factors and stopping at exhaustion}
For the large-formal branch the two-transform return~\cite{Liu,H} has
\[
 M'=M+J-g-g_2+t_2\le M-\sigma,\qquad
 \Delta'=b_2-p_2+a+B_C+\ell_a\ge a+\ell_a.
\]
The main choice is $g=J_++2\sigma$, $\ell_a=\sigma/3$, $a=a_o=0$; fixed amplifier errors use $g=J_++\sigma$, $\ell_a=0$. For a prime-free input no auxiliary-prime average is introduced. First expand the complete last-divisor identity and all exact sixth-power punctures. Only then split the final inducing character according to membership in $\Theta$.

\begin{lemma}\label{lem:exhaustion}
For an actual ordinary child outside $\Theta$, let $D_{\rm ch}$ be its positive-prime affine defect, including the assigned support allowance. If the largest prime-slot width is $\omega$, removal of whole slots until admissibility or exhaustion has squared cost at most
\begin{equation}\label{eq:cap}
 U^{(D_{\rm ch})_+/6+\kappa\omega+\eps}.
\end{equation}
The original central extraction coefficient is preserved.
\end{lemma}
\begin{proof}
Let $z_{\rm rem}$ be the actual remaining width. Each removed slot of width $b_i$ lowers the positive-prime defect by $Kb_i$. Before the last necessary removal their total is less than $(D_{\rm ch})_+/K$; the last slot costs at most $\omega$. If all slots disappear earlier, their total is $z_{\rm rem}$. Thus
\[
 b\le\min\{z_{\rm rem},(D_{\rm ch})_+/K+\omega\},\qquad
 \kappa b\le(D_{\rm ch})_+/6+\kappa\omega.
\]
If a prime remains, the stopping condition is the ordinary moment's length condition. At exhaustion use the prime-free ordinary theorem of~\cite{Liu}, which permits arbitrary bounded ordinary lengths. It is not necessary that $D_{\rm ch}\le0$ in that case. The multiplier already extracted in the central ordinary representation remains outside this estimate; removing a prime does not renormalize either plain factor. The prime-free theorem is a preceding input, so this call is not an induction on the estimate being proved.
\end{proof}

The actual common-support maximum in the return is
\[
 -d-(K_0-k_f)-a_o+J_+
 =\max\{-d-(K_0-k_f)-a_o,-c_C-2a\}\le-\mu_{\rm com}+\delta_1.
\]
Support fixing and last-divisor fixing only decrease the defect. Combining this fact with $\Delta'\ge a+\ell_a$ pays the amplifier's own portion once and leaves
\begin{equation}\label{eq:cost}
 (D_{\rm free}-\mu_{\rm com})_+/6+\eps_{\rm loc}.
\end{equation}
Both Cauchy integrals receive the same paired divisor weight and are bounded separately as in Lemma~\ref{lem:projection}. No bare divisor count is inserted afterwards.

\subsection{Normalization and the remaining branches}
Set $D_0=D_{\rm free}-\mu_{\rm com}$. Direct substitution in~\eqref{eq:completed} gives
\begin{align*}
 R_{\rm loc}&=-f_0-f_1-u-3(T_1+T)-4c-S_1-S-e_1-e\le0,\\
 R_{\rm loc}+D_0/6
 &=-{5f_1\over6}-{19(T_1+T)\over6}-{5v\over6}-4c-S_1-S-e_1-e\\
 &\hspace{7mm}-{\eta_K\over6}-{(K-1)z_C\over6}
                 +{w_f-\mu_{\rm com}\over6}\le0.
\end{align*}
Consequently $R_{\rm loc}+(D_0)_+/6\le0$, which pays~\eqref{eq:cost} for either sign of $D_0$. This includes the cases $z=0$, $w=0$, and complete common support.

For $m_f\le1/100$, use the original Gauss polynomial and $|G(n,h)|\le\sqrt{q_n}$, retaining unit frequencies. The direct calculation in~\cite[the small-formal-length estimate]{H}, after the original coefficient energy and benchmark have been subtracted, gives
\[
 3A_t-\tfrac72+18w+\tfrac7{200}+o(1).
\]
Put $d_0=1-2r_o\in[7/50,9/50]$. Since $A_t\le1-c_Kd_0+W$ and $W<d_0/7$, this is less than
\[
 -\tfrac12+3d_0/K+\tfrac7{200}
 \le-\tfrac{93}{200}+\tfrac{90}{749}<-\tfrac13.
\]
This argument includes negative formal centers and does not apply an ordinary theorem to an unformed Gauss polynomial.

At the first Poisson origin, the character mean vanishes unless the product modulus is powerful. Its divisor-bounded allocations, central square normalization and row mean give $O(U^{1+\eps})$, as in~\cite{Liu}. The second origin has ceiling $A_t-1+4t+o(1)$. The sparse exceptional-character branch has ceiling $A_t-1+3t+o(1)$ after the two matching fixed-character leading means cancel at the exact common scale product. Their respective bounds are
\[
 -c_Kd_0+5W<-(c_K-5/7)d_0<0,\quad
 -c_Kd_0+4W<-(c_K-4/7)d_0<0.
\]
Here $K>7/2$. The exceptional means are handled in this branch, never by Lemma~\ref{lem:exhaustion}. Subunit rectangles retain the original volume alternative.

Integrate in the original weighted Fourier measures and use~\eqref{eq:profile} for tails. This proves the centered bound. To restore $P$, reflect the longer ordinary factor of $P^{\rm cmp}$, retaining its conductor and Euler deletion. Its ordinary total length is at most $3/2-2m$, and
\[
 3/2-2m+K(z+w)\le5/2-4m+Kw<1.
\]
The ordinary input therefore bounds $P^{\rm cmp}E$ on $\mathcal U$. The squared triangle inequality proves the second estimate of~\eqref{eq:norm}, including its cross term. This completes the proof of Theorem~\ref{thm:norm}.

\section{Prime selection and the first mean}\label{sec:selection}
On the main amplitude rectangle
\[
 \mathcal R=[7/20,3/5]\times[2/5,1/2],\qquad g=\delta x,
\]
the detector count has exponent $R_{k_*}(\delta,x)$, and the screened prime inventory has length at least $L_0=21/100$ in $U$ units. The whole-slot selection argument of~\cite[the amplitude-selection lemma]{H} gives residual amplitude
\begin{equation}\label{eq:amplitude}
 a(z)=\delta{xL_0-z/2\over L_0-z},\qquad
 a'(z)={\delta L_0(x-1/2)\over(L_0-z)^2}\le0.
\end{equation}
Choose, at first without strict shortening,
\[
 z_K(m)={1-2m\over K},\quad W_K(m)={c_K(1-2m)\over6},\quad
 W_{\rm new}(m)=W_K(m)+{1\over1000}.
\]
For $41/100\le m\le43/100$,
\begin{align}\label{eq:margins}
 1-2m-7W_{\rm new}&\ge{40523\over6741000}>{3\over500},\\
 1-(5/2-4m+KW_{\rm new})&\ge{15343\over500000}>{3\over100}.\nonumber
\end{align}
Cover this interval by finitely many closed intervals $[m_-,r_o]$ of width less than $1/1000$ and use $W=W_{\rm new}(m_-)$. The first margin remains larger than $1/250$. This supplies legal uniform instances of Theorem~\ref{thm:norm}, rather than combining the maximum width at the left endpoint with the rightmost ordinary length.

\begin{proposition}\label{prop:gain}
On this band, the full first mean gains a net $\nu=1/20000$ in $U$ units over the source amplified-contour envelope, retaining its original physical numerator and both incoming profiles.
\end{proposition}
\begin{proof}
We have $z_K<41/1000$, $W_{\rm new}<1/40$, and $a(z_K)>13/100$. For the last inequality, after multiplication by $L_0-z_K>0$, the difference from $13/100$ is at least
\[
 {21\over10000}-{9z_K\over200}>{51\over200000}>0.
\]
The new first-mean prime gain is $G_{\rm new}=(gz_K+a(z_K)W_{\rm new})/2$. Relative to width $W_K$ it improves by
\[
 G_{\rm new}-G_K={a(z_K)\over2000}>{13\over200000}.
\]
The previous widths are $z_0=2(1-2m)/9$, $W_0=7(1-2m)/54$. Direct substitution in~\eqref{eq:amplitude} shows $G_K\ge G_0$: writing $v=1-2m$, $z_K=sv$, $z_0=s_0v$, $s=1/K$, $s_0=2/9$, one has
\[
 2(G_K-G_0)=g(s-s_0)v+
 {\delta v\over6}\left\{{(1-s)(xL_0-sv/2)\over L_0-sv}
 -{(1-s_0)(xL_0-s_0v/2)\over L_0-s_0v}\right\}\ge0.
\]
Indeed, writing the total saving as $F(s)=G_K$ with $z=sv$, its derivative is
\[
 F'(s)={\delta v\over2}\left\{x-{1\over12}+
 { (1/2-x)L_0(L_0-v)\over6(L_0-sv)^2}\right\}>0,
\]
since $v\le9/50<L_0$ and $x\ge2/5$. This compares against the actual-$m$ source selection, which is at least its uniform contact envelope.
Thus only the extra $1/1000$ is counted as new gain.

Shorten both widths by $1-10^{-4}$. Since $G_{\rm new}<1/100$, this costs less than $10^{-6}$. The inventory remains sufficient because
$2xL_0-z_K>127/1000>W$. Choose whole-slot rounding, support edges, count and norm errors with aggregate first-mean cost below $6\times10^{-6}$, and local height cost below $2\times10^{-6}$. The remaining gain exceeds $5\times10^{-5}$.

For completeness, if $\#\mathscr B\ll U^{R+\eps_R}$ and the selected product has amplitude at least $U^{gz+a(z)w-\eps_A}$, the complete norm gives
\[
 \sum_{u\in\mathscr B}|S_m(u)|^4
 \ll U^{1-2gz-2a(z)w+\eps_N+2\eps_A}.
\]
Hölder gives first-mean loss $(3\eps_R+\eps_N+2\eps_A)/4$. The physical factor $U^{-\theta m}$ stays outside, unchanged. Restriction to $\mathscr B$ occurs after the positive estimate, so the signed coefficient expansion and all cross terms are retained.
\end{proof}

\section{The complete scalar comparison}\label{sec:scalar}
For $k\in[2/3,3/4]$, define
\begin{gather*}
 c={1\over3k},\quad A={5\over6}-\delta,\quad
 D_c=3-(1+2c)x,\quad P_c=2-2(1+c)x+2cx^2,\\
 J_c=AD_c+\delta P_c,\qquad
 R_k=1-\delta+{A\delta P_c\over2J_c}.
\end{gather*}
These are the class-count functions of~\cite{H}. Set $b=b_*$ and
\begin{align}\label{eq:geometry}
 h&={1+b+3\ell\over2},\quad l_y={1+b-\ell\over2},\quad l_x={1-b-\ell\over2},\\
 C_0&=-\tfrac14+b/6+5\ell/4,\quad L=(23\ell-1-b)/18,\quad\theta_0=2\delta x/9,\nonumber\\
 E&=C_0+\delta/2+\delta xL-(3h/4)(1-R_{k_*}),\nonumber\\
 S&=(3R_{k_*}+1)/4+\theta_0/2-17/50.\nonumber
\end{align}
The full unamplified exponent at row scale $U=Z^d$ is $E+(d-h)S$. In particular $E$ already includes the functional-equation and outer normalization factors.

\subsection{Elementary bounds and the algebraic contact}
Differentiation gives
\[
 \partial_kR_k={\delta xA^2(1-x)^2\over3k^2J_c^2}\ge0,\quad
 1-\delta\le R_k\le R_{3/4}\le1-2\delta/3-22\delta x/81.
\]
On $1/8\le\delta\le5/6$, $0\le x\le1/2$,
\begin{equation}\label{eq:derivatives}
 |E_b|<1/5,\qquad 5/16\le E_\ell\le5/4.
\end{equation}
The derivative formulas are $E_b=1/6-\delta x/18-3(1-R)/8$ and $E_\ell=5/4+(23/18)\delta x-9(1-R)/8$.

The polynomial $p$ is strictly decreasing on $[.45,.46]$. With $t=r-.45057$, Taylor expansion gives
\[
 p=.162414095686326-217502.0989101246\,t
                  -81009.52478\,t^2+218582t^3.
\]
At $t=74/10^8$ and $75/10^8$ its signs are positive and negative. Rational substitution in~\eqref{eq:root}--\eqref{eq:b} then yields
\begin{equation}\label{eq:box}
 .45057074<r_*<.45057075,\quad .16879<\ell_*<.16880,
 \quad .09175<b_*<.09185.
\end{equation}
For example, writing $t=r-.45$ gives denominator $7.965-2953.6t+2286t^2$ and numerator $3.375-4058t+3410t^2$ for $\ell_*$. The fixed-$k$ derivatives defining $b_*$ lie respectively in $(21592.47,21592.66)$, $(-4134.33,-4134.28)$ and $(5333.96,5334.03)$. These rational enclosures suffice for every bound below.

Put $a_\Box=\delta(189x-25)/139$ and
\[
 Q=E-{7a_\Box\over108}(2l_y-h).
\]
\begin{lemma}\label{lem:contact}
At $\ell=\ell_*$, $Q\le0$ on $\mathcal R$, with equality only at $(r_*,1/2)$. On this rectangle $S\ge71/300$.
\end{lemma}
\begin{proof}
Let $T=P_cD_c'-P_c'D_c$. Then
\[
 T=4+2c-12cx+2c(1+2c)x^2,\quad (T/D_c^2)'={8(1-c)^2\over D_c^3}>0.
\]
Since $P_c/D_c\ge3/8$, $T/D_c^2<2/3$ and $\delta/A\ge21/29$,
\[
 {T\over[D_c+(\delta/A)P_c]^2}\le{107648\over261075}<5/12.
\]
Thus
\[
 Q_x/\delta=L-{49\over556}(2l_y-h)
       -{3h\over8}{T\over[D_c+(\delta/A)P_c]^2}>0,
\]
using $L>3/20$, $0<2l_y-h<1/8$, $h<4/5$. At $x=1/2$ exact expansion gives
\[
 Q(\delta,1/2)={b_*n_b(\delta,k_*)+\ell_*n_e(\delta,k_*)+n_0(\delta,k_*)
                   \over432(54k_*\delta-75k_*-6\delta+10)}.
\]
Write $D(r)=2286r^2-5011r+1800$ and $K_0(r)=600r^2-1837r+720$. For $\ell(r)$ as in~\eqref{eq:root}, $k(r)=K_0(r)/(3D(r))$ one has
\[
 n_b(r,k(r))=0,\quad
 3D(r)^2\{\ell(r)n_e(r,k(r))+n_0(r,k(r))\}
                   =-1944r(6r-5)^2p(r).
\]
The numerator therefore has a double zero at $r_*$, by~\eqref{eq:b}, and equals $A_Q(\delta-r_*)^2$, where
\[
 A_Q=(-6858k_*+600)b_*+(-702k_*-408)\ell_*+4806k_*-696>2200.
\]
Its denominator is negative on $\mathcal R$. This proves the assertion about $Q$. Finally
$S\ge33/50-3\delta/4+\delta x/9\ge71/300$.
\end{proof}

\subsection{Complementary rectangles and the full box}
Use the rational reference $b_0=.0918$, $\ell_c=.16880$ and count $k=3/4$. Then $h_0=.7991$, $C_{00}=-.0237$, $L_0'=4651/30000$. Direct differentiation gives
\[
 E_{\delta\delta}=-{(3h_0/4)(5/6)^2D_cP_c^2\over J_c^3}<0.
\]
For $\delta\le23/100$ the elementary count bound gives
$E\le-1193/2000000$. For $\delta\ge23/100$, $E$ increases with $x$: indeed
\[
 {T\over[D_c+(\delta/A)P_c]^2}
 \le{856\cdot181^2\over7663^2}<1/2,
\]
whereas $L_0'>3/20$, $h_0<4/5$.
At the two relevant $x$-boundaries the exact expressions are
\begin{align*}
 E(\delta,1/2)&={-2659452\delta^2+2425492\delta-526140
                    \over120000(185-138\delta)},\\
 E(\delta,2/5)&={-201477144\delta^2+167193610\delta-35905500
                    \over600000(2525-1698\delta)}.
\end{align*}
Concavity and the derivative signs at $\delta=7/20,3/5$ imply
\[
 E(\delta,1/2)<-9/50000\ (\delta\le7/20),\qquad
 E(\delta,1/2)<-1/500\ (\delta\ge3/5)
\]
in the range under consideration. To bound the other edge, the numerator of $E(\delta,2/5)+1/1000$, with $y=\delta-2/5$, is
\[
 -201477144y^2+(24965474/5)y-3922976/25
 \le-156000+5000000|y|-200000000y^2<0.
\]
This proves the reference bound $E<-9/50000$ on the three complete rectangles
\[
 [1/8,5/6]\times[0,2/5],\quad [1/8,7/20]\times[2/5,1/2],
 \quad[3/5,5/6]\times[2/5,1/2].
\]
Restoring $b_*$ costs less than $1/100000$ by~\eqref{eq:box}; increasing $\ell$ up to $.16883$ costs at most $3/80000$. Returning to $k_*$ lowers $E$. Hence throughout these rectangles
\begin{equation}\label{eq:complement}
 E<-{53\over400000}.
\end{equation}

We also need a full-box upper bound. Subtract $19/10000$ from $E(\delta,1/2)$ and clear the positive denominator. At $y=\delta-23/50$ its numerator is
\[
 -2659452y^2+(256504/25)y-537677/625
 \le-860+11000|y|-2600000y^2<0.
\]
The preceding monotonicity and small-$\delta$ estimate show $E_{\rm ref}<19/10000$ on the entire box. Restoring the parameters gives
\begin{equation}\label{eq:fullbox}
 E<{779\over400000}<{1\over500}
 \qquad(\ell_*\le\ell\le.16883).
\end{equation}
These arguments cover the boundary edges as well as the interiors.

\subsection{The actual short sums}
For $m\le3/40$, ideal counting in~\eqref{eq:ordinary} gives the original numerator bound $U^{m/2+\eps}$ on each AFE profile. Relative to $N=(3R+1)/4-\theta_0/2$, its first-mean difference is bounded by
\begin{align*}
 \Delta_s&=(R-1)/4+3/80+17\theta_0/40\\
 &\le3/80-\delta/6+43\delta x/1620\le-{1049\over64800}
 \qquad((\delta,x)\in\mathcal R).
\end{align*}
At $d=h$,~\eqref{eq:fullbox} and $h>.79$ give
\begin{equation}\label{eq:short}
 E+h\Delta_s<-{69911\over6480000}<-1/100.
\end{equation}
The actual slope for any $m\le3/40$ is
$R-17/50+\theta_0+(1/2-\theta_0)m$, which is positive on $\mathcal R$ and less than one. Thus the endpoint $d=h$ bounds every licensed smaller $d$. Small rows retain their separate direct estimate. This bound applies to the complete short numerator, with its original outer factors.

\section{Insertion in the physical comparison}\label{sec:global}
For a fixed amplitude class and a single physical dyad, displacement of the contour and a first-mean moment saving $\gamma$ change the exponent by
\begin{equation}\label{eq:physical}
 (\theta-\theta_0)\{d(1-m)-l_y\}-d\gamma,
\end{equation}
as in~\cite[the simultaneous comparison]{H}. Both operations apply to the same sum. The comparison does not assume that the detector and numerator have the same twist or Mellin height.

The AFE is applied at the actual conductor. Writing its index as $a b$, with $a$ supported in $S$ and $q_b\asymp D$, both arms retain the scalar $q_a^{-1/2-\theta}D^{-\theta}$. Their conductor-dependent profiles are handled by the common log-scale Sobolev estimate. On reflection the two AFE arms exchange; the root numbers, conductor powers and Gamma ratios cancel. The entire holomorphic dyad is moved before taking absolute values or conjugating the whole dual arm. Long ordinary factors with $m>1/2$ use the prime-free arbitrary-length theorem, rather than the selected range of Theorem~\ref{thm:norm}.

For comparison, first absorb the independent incoming heights into $V_1,V_2$, then construct the equal-product comparison with those same profiles. An additional Fourier mode acts on the whole product on a given centered side. The fixed-character leading means cancel at the exact common product of scales. Reflection may increase the long factor after a divisor is extracted; its central extraction coefficient must therefore be retained. For the actual top comparison the displayed modulus is its natural radical, and the capacity is the last inequality of~\eqref{eq:domain}.

\begin{proposition}\label{prop:global}
With $b=b_*$ and $\ell=\ell_*+1/40000$, the complete nonprincipal comparison has a scalar reserve greater than $7/10^6$, after a single aggregate global expense smaller than $10^{-7}$.
\end{proposition}
\begin{proof}
From the count bound and~\eqref{eq:derivatives},
\[
 Q_\ell=5/4+(23/18)\delta x-9(1-R)/8+35a_\Box/216
 \le5/4-79\delta/432<593/500.
\]
Lemma~\ref{lem:contact} therefore gives $Q(\ell)<593/20000000$. On $[.41,.43]$, Proposition~\ref{prop:gain} subtracts $h\nu$ in $Z$ units. As $h>.79$, the exponent is less than
\[
 {593\over20000000}-.79\,{1\over20000}=-9.85\times10^{-6}.
\]
For all other ordinary lengths use the original selected and unamplified alternatives of~\cite[the left, central and right envelope comparison]{H}. To make their coverage explicit, put $\lambda=l_y/h\in(.577,.578)$ and $m_0=1-\lambda\in(.422,.423)$. The selected envelope has descending slope at least $7/225$, ascending slope at least $49/432$, and its excess over the contact envelope outside $[.41,.43]$ is at least $7/150000$. The central unamplified range $[3/40,3/8]$ has the same reserve; the outer range $m>9/20$ has the larger excess $7\delta/2160$. These are consequences of the same piecewise-affine envelope, before any class restriction. Its outside-band reserve is thus greater than
\[
 .79\,{7\over150000}-{593\over20000000}-10^{-7}
 ={427\over60000000}>{7\over10^6}.
\]
The complementary classes use~\eqref{eq:complement}; their slope $S$ is positive. The short sums use~\eqref{eq:short}. For $d\le h-1/80$ on $\mathcal R$, Lemma~\ref{lem:contact} and~\eqref{eq:fullbox} give
\[
 E+(d-h)S<{779\over400000}-{71\over24000}
                   =-{1213\over1200000}.
\]
All lie farther from zero than the outside-band reserve.

The low-amplitude, principal elementary and very-small-row branches in~\cite{H} are preserved. Their sufficient scalar conditions follow already from the coarse box $.0917<b<.0919$, $.1687<\ell<.1689$:
\begin{gather*}
 l_x>\ell,\quad l_y-\ell>11b/6,\quad h<2l_y,\quad l_y<1/2,\\
 \ell(1/2+2\ell/h)-11h/54<-1/200,\qquad
 C_0+(1-h)/16<-1/250,\\
 C_0+(1/2+3\ell/2)/50<-1/125,\qquad
 13h/75-l_y/2+1/50<-1/20.
\end{gather*}
They retain the original local corrections and actual support. The affine slopes in the near-$h$ extension are bounded by two, so this extension can be charged inside the aggregate error.

Choose the numerical meshes, strict contour buffers and support edges first. Fix the finite prime list, tests and final finite character group, then every internal Fourier, Sobolev and tail order. If $C_L,C_G$ are the local and still-unpaid global height orders, choose the common external height power below
\[
 \min\left\{{(1/100)\,2\times10^{-6}\over1+C_L},
                       {2\times10^{-8}\over1+C_G}\right\}.
\]
The local term is already paid in $\nu$. All remaining global nonheight errors are chosen together below $8\times10^{-8}$. The total global expense is therefore below $10^{-7}$, as used above.
\end{proof}

\section{Principal recovery and Mellin continuation}\label{sec:continuation}
The preceding scalar reserve concerns the nonprincipal amplitude comparison. The full difference $J_\eta-f_\eta$ also contains the recovery error of the principal signal. We now choose a positive saving common to all its components.

After the numerical mesh has been fixed, choose a finite even prime list with positive lengths $\ell_i$ summing to $\ell$ and satisfying that mesh. Define
\begin{equation}\label{eq:fullsaving}
 r_{\rm np}={7\times10^{-6}\over8},\quad
 r_P=\min\{l_y/20,h/600,(437/1000)\min_i\ell_i\},\quad
 s={\min(r_{\rm np},r_P)\over32}>0.
\end{equation}
The principal contour in~\cite{Liu} gives two bulk savings and a selected-prime recovery saving, respectively
\[
 l_y/20-(20+h)e,\qquad h/600-20e,\qquad\omega_P-20e,
\]
with $\omega_P\ge(437/1000)\min_i\ell_i$. Take $e$ small after the positive $\ell_i$ have been chosen. All three savings are then at least $r_P/2$. In the double residue the factors $Z^{-\ell/6}$ and $Z^{\ell/6}$ cancel before division by the original normalizer $c_SA_T$. Thus the residue gives exactly the same $f_\eta$ as in~\eqref{eq:mellin}.

Choose the normalizing loss and remaining late logarithmic losses below $\min(r_{\rm np},r_P)/64$. This is done by increasing the threshold for the already fixed data. It requires no new prime mesh. After all internal height orders have been fixed, shrink the external height exponent once more so that its complete high-side cost is below $\min(r_{\rm np},r_P)/64$. Finally choose the independent external tail order large enough for the rapidly decreasing Mellin kernels. These operations give
\begin{equation}\label{eq:high}
 |J_\eta(Z)-f_\eta(Z)|\ll_\eta Z^{\beta-c_0-s},\qquad \beta=\beta_*(F).
\end{equation}
The finite prime list and $s$ are common to the target family. Its internal orders, constants and thresholds may depend on $\eta$. In particular,~\eqref{eq:fullsaving} makes no lower bound of $7\times10^{-6}$ for the complete analytic saving. Mesh selection precedes the principal minimum, so the order of choices is noncircular.

\begin{proof}[Proof of Theorem~\ref{thm:main} over $F$]
Suppose $\beta>B$. Choose $0<\eps_{\rm low}<(\beta-B)/4$. The direct input estimate and~\eqref{eq:high} imply
\[
 |J_\eta(Z)|\ll_\eta Z^{B-c_0+\eps_{\rm low}},\qquad
 |f_\eta(Z)|\ll_\eta Z^{\sigma-c_0},\quad
 \sigma=\max\{B+\eps_{\rm low},\beta-s\}<\beta.
\]
The common Mellin identity~\eqref{eq:mellin} therefore continues $H_\eta/L_F^S$ holomorphically to $\Rea s>\sigma$. Since $B>437/500$ and $H_\eta$ is nonzero there, each $L_F^S(s,\eta)$ is zero-free in that common half-plane. The deleted finite Euler factors have no zeros or poles in it. This contradicts the definition of $\beta$ as the supremum over the entire finite-order family. Hence $\beta\le B$.

Only positivity of $s$ is needed: the continuation improves the supposed supremum $\beta$, and does not require $s\ge1/160000$. The larger length $\ell$ determines the direct boundary $B$, while~\eqref{eq:fullsaving} supplies the strict continuation gap.
\end{proof}

\section{Descent to arbitrary number fields}
We finish the field argument, using full local Euler factors. Let $K$ be any number field, $\eta$ a finite-order Hecke character, and $M/K$ its cyclic class field. Put $F=K(\zeta_{12})$ and $E=MF$. The extensions $E/F$ and $E/M$ are abelian. If the inputs of Section~\ref{sec:inputs} hold over $F$, the fixed-field theorem implies
\[
 \zeta_E(s)=\prod_{\psi\in\widehat{\operatorname{Gal}(E/F)}}L_F(s,\psi)
\]
is nonzero in $B<\Rea s<1$. On the other hand the abelian factorization over $M$ is
\[
 \zeta_E(s)=\zeta_M(s)
 \prod_{\substack{\xi\in\widehat{\operatorname{Gal}(E/M)}\\\xi\ne1}}L_M(s,\xi).
\]
The nontrivial factors are holomorphic there, by the standard Hecke theorem and class-field correspondence~\cite{Milne}. Thus a zero of $\zeta_M$ would force a zero of $\zeta_E$. Finally
\[
 \zeta_M(s)=\prod_{\chi\in\widehat{\operatorname{Gal}(M/K)}}L_K(s,\chi)
\]
shows that a zero of $L_K(s,\eta)$ in that strip would force a zero of $\zeta_M$. This is impossible. The line $\Rea s=1$ and the region to its right are covered by the standard nonvanishing theorem and Euler product. The principal pole at $s=1$ creates no zero. This proves the arbitrary-field assertion of Theorem~\ref{thm:main}.

The descent uses the fixed-field result over every $K(\zeta_{12})$; a theorem restricted to the Eisenstein field would not suffice. Its constants need not be uniform in $K$. The numerical boundary is the same because the scalar geometry and the finite-length gain are field independent.

\paragraph{Acknowledgement.}
AI tools assisted in organizing the local derivations, checking algebraic comparisons and preparing the manuscript. The analytic input results are attributed at their points of use.

\end{document}